\documentclass[a4paper,abstracton,11pt]{article}
\usepackage{amsfonts}
\usepackage{amssymb}
\usepackage{amsmath}
\usepackage{amsthm}
\usepackage{xspace}
\usepackage{apacite}

\usepackage{accents}
\usepackage[normalem]{ulem}

\usepackage{booktabs}

\usepackage{graphicx}
\usepackage{xcolor}

\usepackage{subfigure}

\usepackage{a4wide}

\usepackage{setspace}
\newcommand{\X}{{\mathcal{X}}}
\newcommand{\pvar}{{d}}
\newcommand{\IDR}{interpretable optimization rule\xspace}
\newcommand{\IDRs}{interpretable optimization rules\xspace}
\newcommand{\DR}{optimization rule\xspace}
\newcommand{\DRs}{optimization rules\xspace}

\usepackage{authblk}

\allowdisplaybreaks

\begin{document}


\title{A Framework for Inherently Interpretable Optimization Models}

 \author{Marc Goerigk}
 \author{Michael Hartisch\thanks{Corresponding author. Email: michael.hartisch@uni-siegen.de}}



\affil{Network and Data Science Management, University of Siegen,\\Unteres Schlo{\ss} 3, 57072 Siegen, Germany}

\date{}

\maketitle

\begin{abstract}
With dramatic improvements in optimization software, the solution of large-scale problems that seemed intractable decades ago are now a routine task. This puts even more real-world applications into the reach of optimizers. At the same time, solving optimization problems often turns out to be one of the smaller difficulties when putting solutions into practice. One major barrier is that the optimization software can be perceived as a black box, which may produce solutions of high quality, but can create completely different solutions when circumstances change leading to low acceptance of optimized solutions. Such issues of interpretability and explainability have seen significant attention in other areas, such as machine learning, but less so in optimization. In this paper we propose an optimization framework that inherently comes with an easily \IDR, that explains under which circumstances certain solutions are chosen.
Focusing on decision trees to represent \IDRs, we propose integer programming formulations as well as a heuristic method that ensure applicability of our approach even for large-scale problems. Computational experiments using random and real-world data indicate that the costs of inherent interpretability can be very small.
\end{abstract}


\section{Introduction}
\subsection{Motivation}

Due to the existence of explicitly stated models and well defined solution processes, researchers in operations research and optimization can have high confidence in the correctness and usefulness of found solutions. If discrepancies between the real-world application and the found solution occur, the underlying model is usually adapted within a feedback cycle. Overall, this design process gives confidence in the validity of the solution. After obtaining a satisfactory model that reflects all conditions, solutions are expected to be employed in various settings, i.e.  parameters are adapted, yielding new solutions for every special case. 

However, explaining why a solution (e.g. a production plan) of type A is implemented today, while type B was used yesterday is not an easy task. While it is possible to explain the model that is being used, the function that maps a model to a solution is complex and not necessarily interpretable. Furthermore, solvers are able to provide necessary and sufficient optimality conditions, and sensitivity analysis can yield a better understanding of the specifics of an optimal solution. However, most of such checking and validation tools---often communicated in a very formal, mathematical language---are only accessible to experts, that already have a high confidence in the process. Indeed, even understanding a (linear) programming model can quickly become a challenging task, especially for people not familiar with modeling techniques. Users with less mathematical and computer background, e.g. the planner using the optimization software and the workers in charge of implementing the result, may consider the solver a black box. This raises the question whether
a solution is acceptable for the end user, or the object of the optimization process (the worker).
Workers directly affected by the optimized decision may wish to understand them and obtain knowledge about future outcomes in order to asses when one has to work overtime or on weekend shifts. Not being able to answer such questions may lead to discontent and loss of trust and result in weak adoptions of optimized decisions.

We argue that non-theoretical and easily comprehensible explanations have to be available in order motivate actions obtained from optimization programs.  In particular, justifications of \textit{why} a system recommends certain actions are more relevant for non-experts than explanations of \textit{how} the solution is obtained \cite{ji2000use,ye1995impact}. Providing the underlying optimization model might not be useful here, as, ``disclosure alone does not guarantee that anyone is paying attention or is able to accurately interpret the information; more complex information is more likely to be unexamined or misunderstood'' \cite{prat2005wrong}.  
In \cite{smidts2001impact} the authors argue that employees' identification with their company is stronger, the more adequate information they receive about their company and their personal role within the company. Furthermore, effective and transparent communication and explanations play a vital role in fostering  trust and positive employee-organization relationships \cite{rawlins2008measuring,men2014effects,yue2019bridging,doi:10.1287/mnsc.2017.2906}.
Decision models that are complex and difficult to understand can result in skepticism and reluctance to adopt or use the system that was actually designed to support decision making, even if the models have been shown to improve decision-making performance \cite{arnold2006differential,kayande2009incorporating}.
In \cite{kayande2009incorporating} the authors argue that a decision-support system ``must be designed to induce an alignment of a decision maker’s mental model with the decision model''.

Therefore, providing insights into decisions obtained from optimization processes is an important task for managers. Incorporating workers into the feedback cycle of the modeling process can be one important part and providing explanations during this process improves the feedback loop \cite{chakraborti2019plan}. After a model is created workers still are in need of more easily comprehensible explanations of when and why different solutions are considered. Furthermore, end users often have the need of being able to explain decisions obtained from an AI, e.g. for loan \cite{sachan2020explainable,strich2021world} or medical consultants \cite{scott1977explanation,swartout1985explaining,lamy2019explainable}. Even the European ``right to explanation'' tries to ensure that users can ask for an explanation of an algorithmic decision that was made about them \cite{goodman2017european}.

While \emph{explainability} and \emph{interpretability} are often used somewhat synonymously, we briefly discuss commonly accepted differences in their outlook \cite{rudin2019stop}.
A framework is \emph{explainable}, if an easily comprehensible, post-hoc justification of the obtained result is available: When simultaneously looking at the proposed solution and the instance, an  explanation can be given why this specific solution was chosen by the procedure. In contrast, a framework is \emph{interpretable}, if an easily comprehensible rule is available that makes it possible to deduce a result from the instance itself: Because of the characteristics of the instance a (good) solution has to look as follows, i.e.~the mechanics of the procedure can be understood. One of the main problems in these ``definitions'' is the very subjective phrase ``easily comprehensible''. For example linear programming can be viewed as an interpretable framework: given an instance perform the simplex algorithm (the ``easily comprehensible'' rule) and you find the solution. But whether this is in fact easily comprehensible depends on the domain knowledge of the person of interest.

\paragraph{Our Approach.}
In addition to modeling the underlying problem in a suitable manner we impose the task of selecting a well suited \IDR on the modeler: Rather than trying to post-hoc find an explanation of solutions, we want to create \textit{inherently interpretable models} \cite{rudin2019stop} that provide an \IDR, mimicking the optimization process and explaining when to use which solution.
Originating from a set of anticipated scenarios, we provide an \DR that maps any (future) scenario to a solution, giving end users  easy access to the decision making process. In order to ensure comprehensibility we focus on  univariate, binary decision trees of fixed depth, where also the splits at each level of the subtrees are the same. Hence, by performing easy queries on the characteristics of the instance in order to obtain one of the solutions, the optimization process becomes interpretable. Note that we do not claim that the procedure of obtaining the \DR is either explainable or interpretable, but that the \DR provides an interpretable procedure to mimic an optimization process. We show that this places our line of research between so-called $K$-adaptability and decision rule approaches: The former demands the selection of a fixed number of solutions that can be chosen after the disclosure of uncertain parameters, without having a (comprehensible) mapping from scenario to solution. The latter is frequently used in a robust optimization setting, where decision rules approximate the feasible region in order to alleviate the curse of dimensionality, but without focus on comprehensibility.

\subsection{Related Literature}
\paragraph{Explainable and Interpretable AI.}
In AI and optimization the main focus has been the predictive accuracy of models and the quality of found solutions, respectively. However, inexplicable black-box solutions that cannot offer an explanation or justification lack acceptance, 
and in particular for machine learning approaches, the novel field of  explainable AI arose in recent years, as the necessity to comprehend the obtained solution gained importance \cite{paul1993approaches,anderson2018artificial}. For a broad overview we recommend the recent surveys on explainable AI \cite{guidotti2018survey,gilpin2018explaining,vcyras2021argumentative}.

There are various ways of trying to post-hoc explain AI decisions and making them more comprehensible, e.g. counterfactual explanations \cite{martens2014explaining,wachter2017counterfactual}, explaining model predictions by weighting features according to their importance \cite{ribeiro2016should,lundberg2017unified,doi:10.1287/mnsc.2021.4190}, extracting rules explaining the underlying concept \cite{craven1995extracting,martens2007comprehensible}, argumentation \cite{moulin2002explanation,amgoud2009using,atkinson2017towards,vcyras2019argumentation}, or qualitative abstractions \cite{boutilier1994toward,bonet1996arguing}. In \cite{freitas2014comprehensible} the author argues that rather than focussing on accuracy, the comprehensibility of (classification) models should also be considered and in her seminal work \cite{rudin2019stop} Rudin stresses that models should be inherently interpretable rather than post-hoc explainable. 
Several studies exist that research what representation types are well-suited to be comprehensible for human users \cite{miller2019explanation,arzamasov2021comprehensible}. They suggest that decision trees or classification rules \cite{baehrens2010explain,huysmans2011empirical,freitas2014comprehensible} and
linear models \cite{setiono1997neurolinear,parker2015evaluating,ustun2016supersparse}  are quite comprehensible, where sparsity is assumed to be an important factor \cite{gleicher2016framework}. Furthermore, the type of explanation needed highly depends on the targeted audience \cite{arrieta2020explainable}.

\paragraph{Explainable and Interpretable Optimization.}

For problems in the NP complexity class, the feasibility and quality of any solution can be checked in polynomial time. This makes it easy to compare any two solutions, which is known as contrastive explanation \cite{miller2021contrastive}. However, this explanation type demands that an alternative solution has to be present in the mind of the inquirer. In particular for complex domains this is unlikely the case and hence the question of \textit{why} a solution is chosen has to be answered. This question can be answered by operations researchers by pointing at the specific model and the applied algorithm or heuristic, which is very unsatisfying for non-experts.

The possibility of contrastive explanation may have hindered further research into more interpretable models and only few such approaches exist, e.g.  interpretable policies for optimal stopping in stochastic systems \cite{doi:10.1287/mnsc.2020.3592} or a mixed integer programming approach to obtain interpretable tree-based models for clustering \cite{bertsimas2021interpretable}. 
Here, an action (stop or continue) or a class label (cluster membership), respectively, is assigned in an interpretable manner. In contrast, our approach also has to single out \textit{solutions} from an exponentially large solution space, which means that optimization and classification are considered holistically.
In the field of operations research and management science, for scheduling and planning a few recent approaches regarding explainability were published. In \cite{collins2019towards} progress for planning problems towards extracting approximate post-hoc argumentation based explanations is made, where causal relationships are 
extracted and then abstracted. A tool to explain plans to non-technical users is presented in \cite{oren2020argument} that uses formal argumentation and dialogue theory. Explainability in scheduling is discussed in \cite{vcyras2019argumentation}, where a makespan scheduling problem is translated into an abstract argumentation to extract
explanations regarding the feasibility and efficiency. Building upon this approach, a tool that provides interactive explanations in makespan scheduling is presented in  \cite{vcyras2021schedule}. In a more general approach, for solutions obtained via dynamic programming an explaining method was introduces in \cite{erwig2021explainable}, which focuses on the explainability of the program itself. Furthermore, in recent years the importance to provide comprehensible insight gained importance for multi-objective optimization. For a planning environment verbal explanations are generated in \cite{sukkerd2018toward}, that explain the trade-off made to reconcile competing objectives. In a different perspective \cite{corrente2021explainable} use simple decision rules to record the decision makers preference used to iteratively
converge towards the part of the Pareto front containing the best compromise solution, providing insights into the impact of the given answers. Furthermore, an explainable interactive multiobjective optimization method was introduced that supports the decision maker in expressing new preferences in order to improve desired objectives \cite{misitano2022towards}.

\paragraph{Decision Rules and $K$-Adaptability.}
In robust and stochastic multistage optimization, decision variables of later stages---that depend on realizations of uncertain parameters of previous stages---can be approximated using parametric classes of linear or nonlinear functions of the uncertain parameters \cite{georghiou2019decision}. Decision rules specify the reaction to a scenario and therefore can be viewed as an explanation for how a scenario affects a decision.  However, the motivation to use decision rules is not to make decisions of later stages more comprehensible, but to obtain an approximate decision in reasonable time. The selection of a class of functions for the decision rules solely depends on its computational performance and the approximation quality. More recently, predictive prescriptions \cite{bertsimas2020predictive} were introduced that prescribe a decision given the observation of covariate variables, but again not with the focus on obtaining comprehensible prescriptions.

$K$-adaptability allows the decision maker to select $K$ different solutions from which one can be chosen after the disclosure of uncertain parameters. It is used in stochastic \cite{buchheim2019k,malaguti2022k} and robust optimization \cite{bertsimas2010finite,hanasusanto2015k,buchheim2017min} as an approximation scheme with the goal of obtaining good solutions in reasonable time. The focus lies on the selection of the best $K$ solutions, while the decision which solution should be chosen in a specific scenario is postponed to when the scenario is revealed. In particular, no rule is given for the mapping between scenarios and the $K$ solutions. Only in recent works, a connection between $K$-adaptability and decision rules was made \cite{subramanyam2020k,vayanos2020robust}, but not with a focus on obtaining more comprehensible solutions but better computational properties.

\paragraph{Hyper Heuristics.}
The \IDR we want to obtain can be viewed as an easily comprehensible heuristic, which puts our research in the proximity of generation hyper-heuristics \cite{drake2020recent}. Similar to other research fields, the quality of obtained heuristics was mainly in the focus, but already early approaches tried to ensure interpretability by creating representations with limited size \cite{burke2007automatic,nguyen2017pso}. 
Several approaches try to generate a good trade-off between quality and interpretability, e.g.~by post hoc improving interpretability via simplification \cite{nguyen2012computational}, enforcing comprehensible, linear representations \cite{branke2015hyper},
and semantically constraining the arising rule \cite{hunt2016genetic}.
Recent applications of hyper-heuristics where the authors also try to ensure interpretability include routing policies \cite{wang2019novel}, dispatching rules \cite{ferreira2022effective} and online combinatorial optimization problems \cite{zhang2022deep}. All these approaches base their claim of interpretability on a decreased size of obtained rules and expressions.

The most fitting hyper-heuristic classification  \cite{burke2019classification} of the framework proposed in this paper is an offline-learning procedure to obtain a construction heuristic, where in our case ``provision'' would be more fitting than ``construction'', as we immediately provide complete solutions. While most of the approaches to generate heuristics use genetic programming \cite{burke2019classification} or other search algorithms we provide a modeling approach intended as optimization framework. However, as we focus on decision trees as ``candidate heuristic'' \cite{drake2020recent}, we will also compare the optimized decision tree to one obtained via a greedy approach.

\paragraph{Decision Trees.}
In this paper we focus on univariate decision trees in order to provide easily interpretable and comprehensible rules. While finding optimal decision trees is an NP-hard problem \cite{laurent1976constructing}, several heuristics exist that recursively split the feature space in a greedy manner, e.g. in \cite{breiman2017classification,quinlan1986induction}.
Mixed-integer programming formulations \cite{bertsimas2017optimal} can be stated, which can be used to provide good solutions, while they fail to prove optimality for large instances in reasonable time. Other heuristics that are mostly based on pruning include local search \cite{bertsimas2019optimal} and evolutionary algorithms \cite{barros2011survey}. In these methods for obtaining decision trees the classification of the training data is available a-priori. We, on the other hand, want to build an \textit{inherently interpretable model} \cite{rudin2019stop}, i.e.~rather than trying to post-hoc find an explanation, we incorporate the discovery of an explanation into the optimization process by explicitly demanding that a solution needs to have an \IDR. Hence, by predefining the structure of the intended \DR we restrict the solution space instead of dealing with classification errors.

\paragraph{Our Contributions.}
We briefly summarize our main contributions and give an overview of the paper's structure.
\begin{itemize}
\item We present a framework for inherently interpretable optimization with uncertainty in the objective function. This model offers flexibility regarding the \DR which can be chosen to  suitably address the target audience. Additionally, choosing the decision criterion is left to the end user in order to maintain a generally applicable approach (see Section~\ref{sec:model}).

\item We discuss this framework with respect to an \IDR representable via a univariate binary decision tree and focus on the Laplace criterion, which optimizes the set of solutions and the \DR with regard to the average outcome. We examine the computational complexity of the problem and present a greedy approach for such models. A distinct advantage of the greedy method is that in some settings, any solution algorithm for the nominal optimization problem can be utilized (see Section~\ref{sec:discussion}).

\item We provide insights into extensions of the model that can incorporate meta data, uncertainty in the constraints rather than the objective, a multivariate decision tree as \DR, as well as a variant where additional costs are associated with the provision of potential realizations (see Section~\ref{sec:extensions}).

\item In a computational study, utilizing a shortest path setting, we show that the proposed greedy heuristic performs well compared to solving the optimization problem and we provide insights into the quality loss for demanding the existence of an \IDR. In the experiments we utilize both artificial as well as real-world data (see Section~\ref{sec:exp}).
\end{itemize}

\section{An Inherently Interpretable Optimization Model\label{sec:model}}
\subsection{The Framework}

We consider some (deterministic) optimization problem of the form
\[ 
\min_{\pmb{x}\in \mathcal{X}} \pmb{c}^\top \pmb{x}
\]
where $\mathcal{X} \subseteq \mathbb{R}^n$ denotes the set of feasible solutions and $n\in \mathbb{N}$ is the number of variables. Throughout this paper, we use the notation $[n]:=\{1,\ldots,n\}$ to denote index sets. The basic optimization model that we propose is applicable if similar problems are solved repeatedly. We would like to provide an easily comprehensible rule under which circumstances which candidate solution should be used. Finding such an \textit{\IDR} is incorporated into the model rather than post-hoc finding a rule that fits best to the obtained solution. Therefore, the rule's comprehensibility is inherently required, potentially to the cost of no longer being able to obtain globally optimal solutions which might not have an \DR of the demanded form.

So let us assume that the cost vector $\pmb{c}\in \mathbb{R}^n$ is uncertain and a discrete set of candidate cost scenarios $\{\pmb{c}_1,\ldots,\pmb{c}_N\}$ with probabilities $\pmb{p} \in [0,1]^N$, $\sum_{j\in [N]} p_j=1$, is given. Consider the set of all functions $f: \mathbb{R}^n \to \mathcal{X}$ that map cost scenarios to solutions. Let us denote $\mathcal{A}$ as some subset of such functions that we consider as being comprehensible. Our aim is to find an \DR $a \in \mathcal{A}$ that assigns a solution $\pmb{x} \in \mathcal{X}$ to any cost scenario $\pmb{c}\in \mathbb{R}^n$, where in the case of the $j$-th observed cost scenarios we write  $a(j)$ instead of $a(\pmb{c}_j)$ for simplicity. 
The decision-dependent random variable $C(\pmb{x})$ takes a value $C_j(\pmb{x})=\pmb{c}_j^\top\pmb{x}$ with probability $p_j$ for each scenario index $j \in [N]$. Let some decision criterion $\mu_{\pmb{p}}: \mathbb{R}^N \to \mathbb{R}$ be given that aggregates a vector of results to a single objective value (e.g., $\mu_{\pmb{p}}$ can be a risk measure). We can hence evaluate the quality of an 
\DR $a$ by calculating
\[ F_{\pmb{p}}(a) = \mu_{\pmb{p}}(C_1(a(1)),\ldots,C_N(a(N)))\, .\]
The goal is thus to find a rule $a$ that minimizes the function $F_{\pmb{p}}$, i.e., to solve $$\min_{a\in \mathcal{A} }  F_{\pmb{p}}(a)\, .$$
Note that depending on the selected decision criterion, the probability vector $\pmb{p}$ might actually not be needed.

As an example, if the decision criterion $\mu_{\pmb{p}}$ is the variance, function $F_{\pmb{p}}$ is given by 
$$F_{\pmb{p}}(a)=\sum_{j \in [N]}p_j\left(C_j(a(j))-\sum_{k \in [N]}p_kC_k(a(k))\right)^2\, .$$
If we apply the worst-case criterion, $F_{\pmb{p}}$ is given by
$$F_{\pmb{p}}(a)=\max_{j\in [N]} C_j(a(j))\, ,$$
while in the following we often focus on the Laplace criterion, i.e. $F_{\pmb{p}}(a)=\mathcal{L}(a)=\sum_{j\in [N]}C_j(a(j))$ where we assume that all scenarios have the same probability (note that minimizing the sum of values is equivalent to minimizing the expected value in this case).

As the main goal is to obtain an \IDR, the structure of allowed rules in $\mathcal{A}$ is crucial. 
One key aspect in this paper is to restrict the image of the \DR, i.e.~$\{a(j) \mid j \in [N]\}\subseteq\X$, to contain only a small number of potential solutions. Furthermore, the mapping of a scenario to a solution has to be easily comprehensible, given by a straightforward rule. 

This marks the central difference between our approach and $K$-adaptability as in \cite{buchheim2019k}. While the $K$-adaptability approach also restricts the cardinality of the image of the \DR (hence the $K$ in the name), it does not put any further constraints on it to ensure comprehensibility, i.e., any mapping from $[N]$ to $[K]$ is allowed. This means that usually, the best out of the $K$ candidate solutions is assigned to each scenario. In the following, we refer to the $K$-adaptability approach with the Laplace criterion $\mathcal{L}(a)$ as the min-sum-min approach, in analogy to \cite{buchheim2017min}.

For calculating an \DR, some data set of scenarios $\{\pmb{c}_1,\ldots,\pmb{c}_N\}$ is required for planning. In practice, it may happen that different---not yet observed---scenarios are encountered, and we still expect the rule to perform well. For this reason, we usually distinguish between training data (i.e., the $N$ scenarios available for optimization) and test data (i.e., a set of scenarios not known during optimization). In particular our approach is suitable whenever the true scenario is disclosed and \textit{observable} before a decision has to be made and the used solution is desired to be obtainable via an easily comprehensible rule. Finding such an \DR is therefore a long-term planning decision, where the \DR should be applicable in several future scenarios. Application examples include  routing problems, disaster management and workforce scheduling.

\subsection{Illustrative Example}
Consider a type of selection problem, where $n=5$ projects are given and we would like to create a portfolio consisting of exactly $p=2$ of these projects. More formally, we have $\X = \{ \pmb{x}\in\{0,1\}^n : \sum_{i\in[n]} x_i = p\}$. Each project $i\in[n]$ has an associated cost, which we would like to minimize. However, these costs vary and depending on their true realization we want to select different projects. To this end we want to provide an \DR that allocates solutions to occurring cost scenarios. As we aim for an easily \IDR we use a univariate, binary tree of fixed depth with the same splits in each level of the tree. By fixing the depth of the tree to two, at most four different solutions $\pmb{x}^1,\pmb{x}^2,\pmb{x}^3,\pmb{x}^4\in \X$ can be attained at the leaves. We assume there are $10$ cost scenarios with equal probability (first columns of Table~\ref{tab:ex2}), and we would like to minimize the expected costs.

\begin{figure}[htb]
\begin{center}
\subfigure[Splits in $(c_2,c_3)$ plane.]{\includegraphics[width=0.4\textwidth]{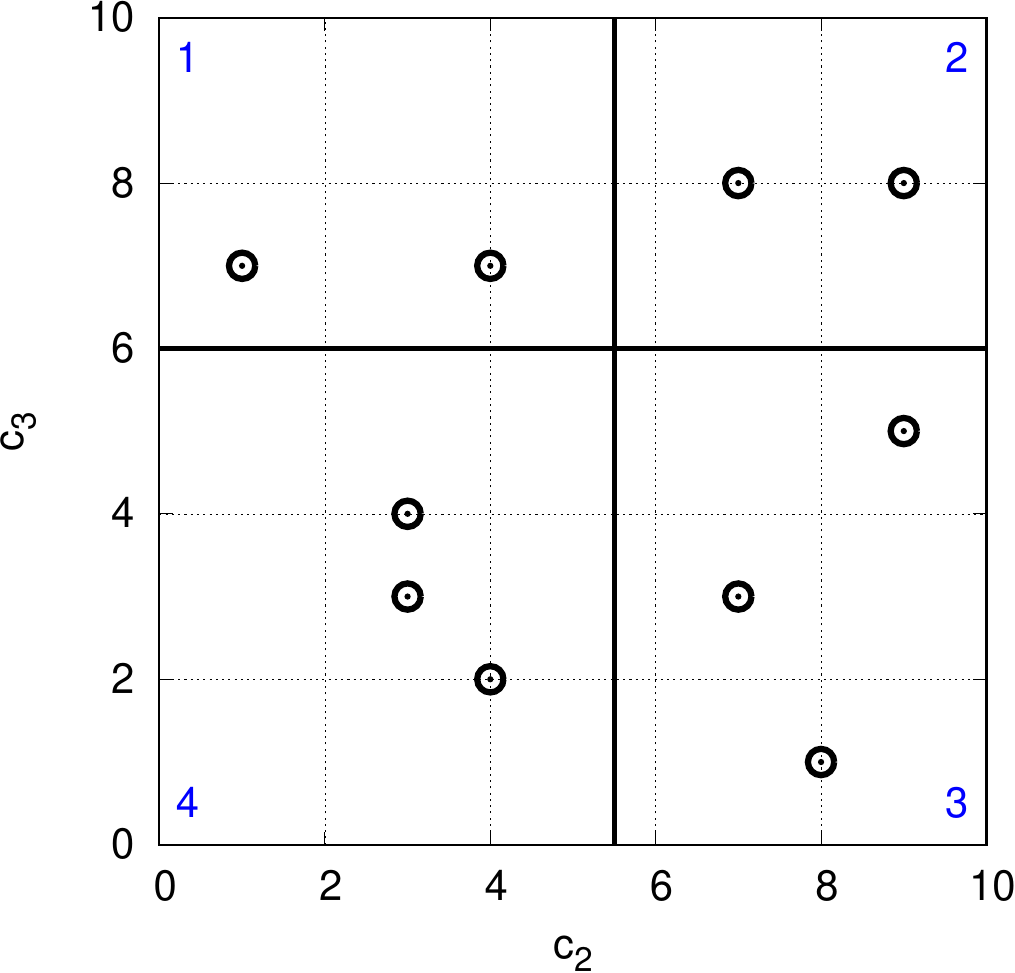}}
\hfill
\subfigure[Decision tree.]{\includegraphics[width=0.5\textwidth]{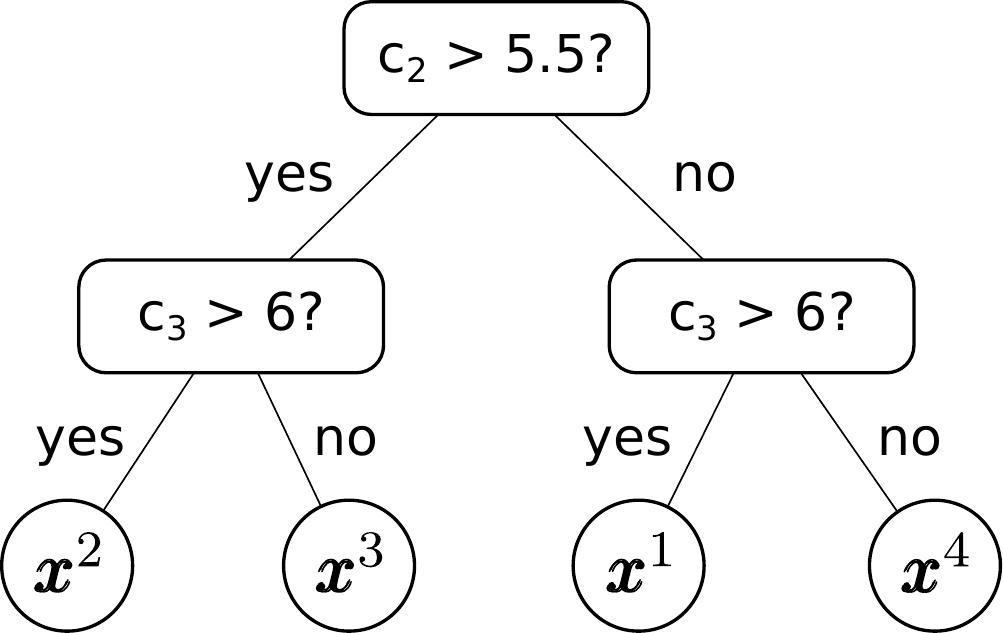}}
\caption{Example \DR.}\label{fig:ex}
\end{center}
\end{figure}

As the characteristics of the desired \DR are already established, this decision tree has to be filled with life: Good splits as well as the proposed solutions associated to each leaf need to be found. \emph{This task is performed via our optimization framework}, which yields the two splits $c_2 > 5.5$ and $c_3 > 6$, as well as the four solutions $\pmb{x}^1=\{2,5\}$, $\pmb{x}^2=\{1,5\}$, $\pmb{x}^3=\{3,5\}$, and $\pmb{x}^4=\{2,3\}$. Each scenario belongs to exactly one category defined by the selected splits, see Figure~\ref{fig:ex}, where we show the projections of the ten given scenarios onto the $(c_2,c_3)$ plane. We associate a solution $\pmb{x}^k$ with each category, see the numbers in the four corners of the plot. For example, if $c_2 > 5.5$ and $c_3 > 6$, we make use of solution $\pmb{x}^2$. These are rules that are easy to communicate and to check without further calculations.

The corresponding costs are given in Table~\ref{tab:ex2}, along with the allocation choice made by the \DR. Note that in scenarios 6 and 8 none of the four allocated solutions is optimal.

\begin{table}[htb]
\caption{Project costs, solution costs and choice.\label{tab:ex2}}
\begin{center}
\begin{tabular}{r|rrrrr||rrrr|r|r|r}
& \multicolumn{5}{c||}{Project Costs}& \multicolumn{4}{|c|}{Costs of Solution} & Allocated & Optimal& Cost of Inter-\\
Scenario& $c_1$ & $c_2$ & $c_3$ & $c_4$ & $c_5$& $\pmb{x}^1$ & $\pmb{x}^2$ & $\pmb{x}^3$ & $\pmb{x}^4$ &  Solution & Solution& pretability\\
\hline
1 & 4 & 7 & 8 & 6 & 4 & 11 & 8 & 12 & 15 & $\pmb{x}^2=\{1,5\}$&$\{1,5\}$&0\\
2 & 6 & 7 & 3 & 10 & 2& 9 & 8 & 5 & 10 &   $\pmb{x}^3=\{3,5\}$& $\{3,5\}$&0\\
3 & 3 & 8 & 1 & 10 & 4& 12 & 7 & 5 & 9 &   $\pmb{x}^3=\{3,5\}$& $\{1,3\}$&1\\
4 & 7 & 1 & 7 & 3 & 7 & 8 & 14 & 14 & 8 &  $\pmb{x}^1=\{2,5\}$& $\{2,4\}$&4\\
5 & 2 & 9 & 8 & 10 & 3& 12 & 5 & 11 & 17 &   $\pmb{x}^2=\{1,5\}$& $\{1,5\}$&0\\
6 & 1 & 3 & 4 & 8 & 6& 9 & 7 & 10 & 7 &  $\pmb{x}^4=\{2,3\}$& $\{1,2\}$&3\\
7 & 10 & 3 & 3 & 9 & 10& 13 & 20 & 13 & 6 & $\pmb{x}^4=\{2,3\}$& $\{2,3\}$&0\\
8 & 8 & 4 & 7 & 2 & 3 & 7 & 11 & 10 & 11 &$\pmb{x}^1=\{2,5\}$& $\{4,5\}$&2\\
9 & 10 & 4 & 2 & 10 & 5& 9 & 15 & 7 & 6 &  $\pmb{x}^4=\{2,3\}$& $\{2,3\}$&0\\
10 & 8 & 9 & 5 & 6 & 1 & 10 & 9 & 6 & 14 &  $\pmb{x}^3=\{3,5\}$ & $\{3,5\}$&0
\end{tabular}
\end{center}
\end{table}

This \IDR that only queries the true cost of projects 2 and 3  is in this example even able to allocate the best of the four available solutions in each scenario but not necessarily the global optimal solution. This is not necessarily always the case. In comparison to a min-sum-min solution, there are two ways how our approach may result in higher costs. On the one hand, the rule may not be able to allocate optimal solutions to every scenario due to its simple (but comprehensible) structure. On the other hand, even if we allocate the best of the $K$ solutions to each scenario, the choice of solutions may still result in an increase in objective value compared to selecting an individual solution in every scenario. Furthermore, note that this only holds when we consider the training data. On additional testing data, which was not used during the optimization, the performance of a solution obtained from an interpretable model may even be better than that of a min-sum-min solution.

\section{Model Discussion}
\label{sec:discussion}

\subsection{Mixed-Integer Programming Formulation for Decision Trees}
\label{sec:formulation}

We first introduce a mixed-integer programming formulation for the problem of minimizing $F_{\pmb{p}}(a)$ when we restrict the \DR to the form of a univariate, binary decision tree of fixed depth. The application of such a tree can be seen as a series of yes/no queries, to reach a leaf that corresponds to a solution. Let $Q\in \mathbb{N}$ be the number of queries, called \textit{splits}, and thus $K=2^Q$ is the number of leaves of the decision tree.  The $k$-th leaf is associated with solution $\pmb{x}^k\in\mathcal{X}$. Note that we do not demand that these solutions have to be different. This might result in an implicit downsizing of the tree in case these are neighboring leaves as then the final split (of this branch) can be left out as it results in the same solution anyway.

To ensure comprehensible rules, we let the split at each level $q$ of the decision tree be the same in each subtree. Additionally, only a single entry of the cost scenario vector $\pmb{c}_j \in \mathbb{R}^n$ of the materialized scenario $j\in [N]$ is queried and compared to a threshold value $b^q\in \mathbb{R}$.  Binary variables $\pvar_i^q$ indicate whether entry $i$ of the cost vector is queried in split $q$ of the decision tree. Let $a_j^k$ be the indicator whether in scenario $j\in[N]$ the traversal of the decision tree reaches leaf $k \in [K]$.

Each leaf can be described via a binary representation of the split outcomes. To this end let $S_q\subseteq[K]$ contain the leaves $k$ where in binary representation of $k-1$, the $q$th digit is 1. For example with $Q=2$ splits and $K=4$ leaves, we have $S_1 = \{2, 3\}$ and $S_2 = \{1, 3\}$ using the binary representation as shown in Table~\ref{Tab::BinaryRep}.
\begin{table}[h!]
 \caption{Binary representations \label{Tab::BinaryRep}}
\centering
\begin{tabular}{r|cc|l}
$k-1$ & $q=1$ & $q=2$ & contained in \\
 \hline
 0 & 0 & 0 & -\\
 1 & 0 & 1 & $S_2$\\
 2 & 1 & 0 & $S_1$\\
 3 & 1 & 1 & $S_1$, $S_2$
 \end{tabular}
\end{table}

Using this notation, we propose the following mixed-integer formulation, that takes the cost scenarios $\{\pmb{c}_1,\ldots,\pmb{c}_N\}$ as input and yields an \DR of the demanded structure (decision tree of depth $Q$, univariate splits, same splits on level $q$ in each subtree) by \textit{simultaneously} selecting solutions ($\pmb{x}^k$) for the leaves and splits ($d_{i}^q$, $b^q$) at each level.
\begin{subequations}
\label{Eq::MainModel}
\begin{align}
\min \ & F_{\pmb{p}}(a) = \mu_{\pmb{p}}\left(\sum_{k\in[K]} a^k_1 \pmb{c}_1^\top \pmb{x}^k, \ldots, \sum_{k\in[K]} a^k_N \pmb{c}_N^\top \pmb{x}^k\right) \label{Eq::MainModel::Obj}\\
\textnormal{s.t.}\ 
&\sum_{k \in [K]} a^k_j=1&& \forall j \in [N] \label{Eq::MainModel::OneSolutionPerScenario}\\
&\sum_{i\in [n]} \pvar_i^q = 1 &&\forall q \in[Q] \label{Eq::MainModel::Explainable1}\\
&  \sum_{i\in[n]} c_{ji} \pvar_i^q \le b^q 
+  M\sum_{k\in S_q}a^k_j && \forall q \in[Q],\ j\in[N] \label{Eq::MainModel::Threshold1}\\
& b^q+\epsilon -M(1-\sum_{k\in S_q} a^k_j)  \le \sum_{i\in[n]} c_{ji} \pvar_i^q && \forall q\in[Q],\ j\in[N] \label{Eq::MainModel::Threshold2}\\
& \pmb{x}^k \in \mathcal{X}&&\forall k \in [K] \label{Eq::MainModel::Domain}\\
&b^q\in [\min_{j\in[N], i\in[n]} c_{ji}-\epsilon,\max_{j\in[N], i\in[n]} c_{ji}]\\
&\pvar_i^q \in \{0,1\}&&\forall i \in [n],\ q\in[Q] \label{Eq::MainModel::Explainable2}\\
&a^k_j \in \{0,1\}&&\forall k \in [K],\ j \in[N]
\end{align}
\end{subequations}
While Constraints \eqref{Eq::MainModel::OneSolutionPerScenario} demand that  each scenario is mapped to exactly one solution, Constraints \eqref{Eq::MainModel::Explainable1} ensure comprehensibility by restricting each split  in the decision tree to query only a single parameter. Constraints \eqref{Eq::MainModel::Threshold1} and \eqref{Eq::MainModel::Threshold2} describe each split of the decision tree and mirror the traversal of a scenario through the different levels: Assume $b^q$ and $\pmb{\pvar}^q$ are already fixed. Then for each scenario $j \in [N]$ and split $q\in[Q]$ either $\sum_{i \in [n]}c_{ji}\pvar_i^q \leq b^q$ or $\sum_{i \in [n]}c_{ji} \pvar_i^q > b^q$ holds. Depending on the outcome the decision tree is traversed to the right or to the left, and hence either one of the leaves to the right or to the left of the current node has to be reached which is indicated via $\sum_{k\in S_q} a^k_j$ for sufficiently large $M$ and sufficiently small $\epsilon$. In particular, with $\epsilon \leq \frac{1}{2} \min_{j,k \in [N],j\neq k, i \in[n]} |c_{ji}-c_{ki}|$ we can ensure the existence of reasonable splits, as splits only have to distinguish between the costs of a single item. Furthermore, as we restrict the thresholds values $b^q$ to the relevant region $[\min_{j\in[N],i\in[n]} c_{ji}-\epsilon,\max_{j\in[N],i\in[n]} c_{ji}]$, using $M \geq 2\max_{j\in[N],i\in[n]} \vert c_{ji}\vert+\epsilon$ suffices in order to ensure the operating principle of Constraints~\eqref{Eq::MainModel::Threshold1} and \eqref{Eq::MainModel::Threshold2}.

Note that this model chooses an arbitrary optimal solution $\pmb{x}^k$ for each leaf, if more than one optimal solution exists. For the sake of communicating resulting decision trees, it may be beneficial to check for each $\pmb{x}^k$ selected this way if it is also optimal in other leaves of the tree in a post-processing step, or to generate a pool of candidate solutions instead of only a single optimal solution (see, e.g., \cite{tsai2008finding}).

The optimization criterion \eqref{Eq::MainModel::Obj} as well as the optimization domain \eqref{Eq::MainModel::Domain} remain general and the type (linear, quadratic, etc.) and thus the tractability of this model vary from case to case.

As a special case, we consider the Laplace criterion, i.e., we are interested in \DRs yielding solutions that perform well on average. 
This optimization criterion can be realized within the scope of Model \eqref{Eq::MainModel} by using the quadratic objective function $\mathcal{L}(a)=\sum_{k \in [K]}\sum_{j \in [N]} a^k_j\sum_{i \in [n]}c_{ji} x_i^k$ that adds up the costs of the solutions chosen in every scenario. If a linear model is of interest, linearization options are possible, e.g.\ by introducing auxiliary variables $z_j \in \mathbb{R}_{\geq 0}$, using objective function $F_{\pmb{p}}(a)=\sum_{j \in [N]} z_j$, and adding constraints
\begin{equation}
\sum_{i \in [n]}c_{ji} x_i^k\leq z_j + M'(1-a_j^k) \qquad \forall k \in [K],\ j \in [N]
\end{equation}
which force $z_j$ to take the cost of the unique solution $k \in [K]$ associated with scenario $j \in [N]$ according to the decision tree. Here, the constant $M'$ needs to be a problem-specific upper bound on the objective values of optimal solutions in each scenario.

\subsection{Notes on Complexity}

In the proposed model, Restrictions~\eqref{Eq::MainModel::Explainable1} ensure that each split within the decision tree only considers a single dimension of the cost vector $\pmb{c}$. While this feature is used for the sake of comprehensibility, it also gives a computational benefit. Given a constant tree depth $Q$, there are $O(n^Q)$ possibilities to set variables $\pvar^q_i$. Once we fix which dimension $i$ to query in split $q$, we can enumerate all possible values for $b^q$. Let us assume that $d^q_i = 1$ for some fixed $q$ and $i$. Consider the set $\{ c_{ji} : j\in [N] \}$ of values in dimension $i$. Let us sort the values of this set, and consider all midpoints between two subsequent sorted values. Then, this set of midpoints contains an optimal choice for $b^q$, as we can enumerate all possible separations of values in this dimension. Hence, there are $O(N)$ possibilities to set $b^q$.

If all $\pvar^q_i$ and $b^q$ variables are fixed, variables $a^k_j$ become fixed as well. Hence, the remaining optimization problem is to minimize the criterion $\mu_{\pmb{p}}(\sum_{k\in[K]} a^k_1 \pmb{c}_1^\top \pmb{x}^k, \ldots, \sum_{k\in[K]} a^k_N \pmb{c}_N^\top \pmb{x}^k)$ with a fixed assignment of solutions to scenarios. We refer to these problems as \textit{subproblems}. If subproblems can be solved in $O(T)$, this means that the interpretable optimization problem can be solved in $O(n^QN^QT)$.

In the special case of the Laplace decision criterion, we have that 
\begin{align*}
\mu_{\pmb{p}}(\sum_{k\in[K]} a^k_1 \pmb{c}_1^\top \pmb{x}^k, \ldots, \sum_{k\in[K]} a^k_N \pmb{c}_N^\top \pmb{x}^k)
&= \sum_{j\in[N]} \sum_{k\in[K]} a^k_j \pmb{c}_j^\top \pmb{x}^k \\
&= \sum_{k\in[K]} \sum_{j\in[N]} a^k_j \pmb{c}_j^\top \pmb{x}^k \\
&= \sum_{k\in[K]} \left( \sum_{i\in[n]} \left( \sum_{j\in[N] : a^k_j=1}  c_{ji} \right) x^k_i \right)\, .
\end{align*}
This means that we need to solve $K$ separate nominal problems, where the costs in each such problem are given by the sum of costs of scenarios that have been assigned to this $k\in[K]$. Hence, the interpretable optimization problem can be decomposed to $O(n^QN^QK)$ nominal problems. In general, if $\mu$ is decomposable (see \cite{jesus2014survey}), it is possible to consider $K$ separate problems to solve the subproblem.

\subsection{Greedy Heuristic}
\label{sec:greedy}

The connection to decision trees as they are used in classification problems suggests the use of a greedy heuristic to solve the interpretable optimization problem. The main idea is that instead of considering all $n^QN^Q$ possible combinations to determine splits in the decision tree, we determine the splits one level after the other.

We first enumerate all possible splits on the first level. That is, we solve a problem with $Q=1$ by enumerating all possible values for variables $\pvar^1_i$ and for the threshold $b^1$. We hence solve $O(nN)$ subproblems in this first step.

We then fix the best of these candidate splits and consider the next level. That is, we solve a problem with $Q=2$, where we only need to enumerate the possible values for variable $\pvar^2_i$ and $b^2$. Again, there are $O(nN)$ subproblems to be solved.

By repeating this procedure until a desired depth $Q$ is reached, a feasible solution has been found in $O(QnNT)$. Note that this means the depth $Q$ has been moved from being an exponent in the full enumeration approach, to being a factor in the heuristic. Further note that by construction, the greedy heuristic is optimal for $Q=1$.

In case of the Laplace criterion, the first subproblem consists of solving two nominal problems, the second subproblem consists of solving four nominal subproblems, and so on. In total, we need to solve $O(nN\sum_{q=1}^Q 2^q) = O(nNK)$ nominal problems to construct a greedy decision tree.

Of course other heuristic approaches, e.g.~evolutionary algorithms as frequently used in generation hyper-heuristics \cite{drake2020recent}, are conceivable to obtain a decision tree of the demanded structure.

\section{Model Extensions}
\label{sec:extensions}

In this section, we discuss possible extensions of the basic interpretable optimization formulation presented in Section~\ref{sec:formulation}.

\subsection{Using Meta Data}
In order to obtain \DRs that are applicable in practice, the input needed to evaluate the rule and retrieve the suggested solution has to be available. So far, the cost value in each dimension of the realized scenario functions as the input for the rule. The ability to retrieve this information at any point in time can be questionable and highly depends on the use case itself. But even if the cost parameters cannot be observed, scenarios---in particular in the case of real world data---often contain more information than seemingly useful for the problem at hand. For example traffic data might not only contain travel times, but also information regarding time and date and even could be augmented to contain weather data or other information that potentially has an impact on the decision: In a shortest path setting it is intuitively clear that an optimal path from the eastern to the western part of a town can depend on the daytime, the day of the week or the weather conditions. Hence, when trying to find easily comprehensible, high quality rules for future scenarios this kind of information should not be dismissed.

In order to utilize meta data, the  general model proposed in Section~\ref{sec:discussion} can be used by expanding the problem domain $\mathcal{X}$ by $m$ dummy variables, where $m \in \mathbb{N}$ is the dimension of the available meta data. Each dummy variable $x_{n+i}$ for $i\in[m]$  should be forced to zero by demanding $x_{n+i}=0$ in $\mathcal{X}$ in order to stop these non-cost information to affect the objective function. The entry $c_{j,n+i}$ for meta information $i$ in scenario $j$ does not contain the ``cost'', but the realized value in this observed scenario, e.g. if for each scenario the day of the week is known, then $c_{j,n+i}$ is a value between 1 and 7.

While implementing this augmentation can be done in a straight-forward manner, expert knowledge of the underlying problem is useful in order to obtain concise \DR. Picking up on the shortest path setting, one useful split in the decision tree could try to separate rush hour times. Let us assume that rush hour is from 6am to 9am and the time stamp of each scenario is stored in a minute-by-minute precision starting at midnight represented by the value 0 until one minute before midnight represented by the value 1439. In order to determine whether a scenario occurs within the rush hour, two decision tree levels are needed, where the first one checks whether it is after 6am, i.e. $c_{j,n+i} \geq 360$ and a second one checks $c_{j,n+i} \leq 540$. However, if the relevance and characteristics of the rush hour are known, the separation of such scenarios could be achieved by a single split at the cost of having to pre-process the data. In the described case, one option could be to shift the time-stamp interpretation such that time stamp 0 represents scenarios at  6am, making the split $c_{j,n+i} \leq 180$ sufficient in order to identify a rush hour scenario. In fact, it might even be replaced by a 0/1 attribute to indicate if this scenario takes place within the rush hour or not.  Another way could be to partition the time stamp information into eight clusters with each three hour duration. This way, instead of having a single data point ``time stamp'', eight indicator values can be used to represent the clusters. Thus, 6am to 9am scenarios could be determined by splitting on the third entry, only. 

\subsection{Uncertainty in Constraints}
So far the changing parameters only affected the quality and not the feasibility of a solution, i.e.\ the feasible set of solutions $\mathcal{X}$ remained the same in all scenarios. If a scenario also affects the feasible region, each scenario $j$ induces a different feasibility set $\mathcal{X}^j$, where the $k$-th solution $\pmb{x}^k$ only has to satisfy the feasibility sets of the scenarios $j$ it is assigned to via $a^k_j$. We denote 
$$
\mathcal{X}^{k,j}(\pmb{a})=
\begin{cases}
\mathbb{R}^n &\textnormal{, if } a^k_j=0\\
\mathcal{X}^j & \textnormal{, if } a^k_j=1
\end{cases}
$$
and demand that $\pmb{x}^k \in \mathcal{X}^{k,j}(\pmb{a})$. Assuming that $\mathcal{X}^j$ can be stated via a linear constraint system $A^j \pmb{x}\leq \pmb{b}^j$, $\pmb{x}^k \in \mathcal{X}^{k,j}(\pmb{a})$ can be enforced via
$$A^j \pmb{x}^k\leq \pmb{b}^j + (1-a^k_{j})\pmb{M}\quad \forall k\in[K],\ j\in[N]$$
with $\pmb{M}$ being a vector with sufficiently large entries, such that in case of $a^k_j=0$ these constraints are trivially fulfilled. Note, however, that the existence of appropriate values for $\pmb{M}$ depends on the structure of $\mathcal{X}$.

If with probability $\epsilon \in [0,1]$  a solution is allowed to not obey the respective feasibility set, additional auxiliary variables $g_{j} \in \{0,1\}$ can be used with $\sum_{j \in [N]} p_j g_j \leq \epsilon$ and the feasible region is given by
$$A^j \pmb{x}^k\leq \pmb{b}^j + (1-a^k_j+g_j)\pmb{M} \quad \forall k\in[K],\ j\in[N]\, .$$

\subsection{Multivariate Decision Trees}

With the goal of obtaining easily \IDRs we restricted ourselves to univariate decision trees (cf. Constraints \eqref{Eq::MainModel::Explainable1}  and \eqref{Eq::MainModel::Explainable2}). Depending on the use case, multivariate decision trees might also be appropriate. A model that can be used to obtain a binary decision tree that allows more general linear splits can be formulated using the following constraints to replace (\ref{Eq::MainModel::Explainable1}-\ref{Eq::MainModel::Threshold2}).
\begin{subequations}
\label{Eq::MultivariateModel}
\begin{align}
&\sum_{i\in [n]} \pvar_i^q \leq P &&\forall q \in[Q] \label{Eq::MultivariateModel::RestrictNonZeros}\\
& \vert \ell_i^q \vert \leq M_\ell d_i^q &&\forall q \in[Q] \label{Eq::MultivariateModel::RestrictNonZeros2}\\
&\sum_{i\in [n]} \vert \ell_i^q \vert =1 &&\forall q \in[Q] \label{Eq::MultivariateModel::RestrictNonZeros3}\\
&  \sum_{i\in[n]} c_{ji} \ell_i^q \le b^q 
+  M\sum_{k\in S_q} a^k_j && \forall q \in[Q],\ j\in[N] \label{Eq::MultivariateModel::Threshold1}\\
& b^q+\epsilon -M(1-\sum_{k\in S_q} a^k_j)  \le \sum_{i\in[n]} c_{ji} \ell_i^q && \forall q\in[Q],\ j\in[N] \label{Eq::MultivariateModel::Threshold2}\\
&\ell_i^q \in \mathbb{R}&&\forall i \in [n],\ q\in[Q] 
\end{align}
\end{subequations}

In order to ensure comprehensibility of the resulting decision tree we uphold in \eqref{Eq::MultivariateModel::RestrictNonZeros} and \eqref{Eq::MultivariateModel::RestrictNonZeros2} that at most $P \in \mathbb{N}$ non-zero entries are allowed in each $\pmb{\ell}^q$, i.e., we still aim for sparse rules. Furthermore, we use a normalization constraint $\sum_{i\in [n]} \vert \ell_i^q \vert =1$ to mitigate numerical issue arising from the use of $\epsilon$ and $M$ in the same constraints. Only allowing integer coefficients can further reduce numerical issues, to the expense of the model's computability.

When using the enumeration approach or the greedy heuristic for multivariate decision trees, more combinations need to be taken into account. If $P$ is a constant, we can still enumerate all possibilities for binary variables $d^q_i$, of which at most $P$ can be active. Values for $\ell^q_i$ and $b^q$ can then be enumerated by considering all $P$-dimensional hyperplanes defined through the costs $c_{ji}$. However, such approaches become computationally too expensive with higher values of $P$.

\subsection{Decision Trees with Varying Splitting Criteria depending on the Subtrees} 
So far the splits used in the decision tree were the same at every level of the tree. This was done in order to maintain the highest possible level of comprehensibility as in this setting the same $Q$ split criteria have to be checked regardless of the scenario. However, having different splits depending on the current region reached while traversing the decision tree can also be of interest, as solutions with higher quality are obtainable with only a slight decrease of comprehensibility. Having the number $K$ of considered solutions  and thus the number of leaves of the decision tree fixed, we now have to decide on the splitting criterion in each of the $K-1$ internal nodes. In order for one leaf to be reached, the splits of the $Q$ nodes traversed to reach this leaf have to be set accordingly.

Following the notation as in \cite{bertsimas2017optimal} we denote as $A_L(k)$ the set of ancestors of the leaf representing solution $k$ whose left branch has to be chosen in order to reach this leaf, and $A_R(k)$ are the remaining nodes on the path from the root to this leaf. By keeping the notation that $a_j^k$ indicates whether solution $k$ is chosen in scenario $j$, but changing $q$ to the index of an internal node rather than the level of the splitting node, Constraints \eqref{Eq::MainModel::Threshold1} and \eqref{Eq::MainModel::Threshold2} (or for the multivariate case Constraints \eqref{Eq::MultivariateModel::Threshold1} and \eqref{Eq::MultivariateModel::Threshold2}) have to be replaced by the following two constraints
\begin{subequations}
\begin{align}
&  \sum_{i\in[n]} c^j_i p_i^q \le b^q 
+  M(1-a^k_j) && \forall k \in[K],\ q \in A_L(k),\ j\in[N] \\
& b^q+\epsilon  \le \sum_{i\in[n]} c^j_i p_i^q+  M(1-a^k_j) && \forall k \in[K],\ q \in A_R(k), \ j\in[N]\, .
\end{align}
\end{subequations}
With Constraint \eqref{Eq::MainModel::OneSolutionPerScenario} each scenario is associated with exactly one leaf and with the above constraints the splitting criteria along the path leading from the root to this leaf have to be met.

Using such decision trees in the enumeration method or the greedy heuristic is directly possible, but means that more variable combinations have to be enumerated. In case of the enumeration method, $O(Q^Qn^QN^Q)$ subproblems need to be solved, while the greedy heuristic needs to consider $O(KnN)$ subproblems.

\subsection{Preparation Costs}
In the proposed model it is assumed that the costs of implementing solution $k$ in scenario $j$ are represented via $\pmb{c}^\top_j\pmb{x}^k$. However, providing the option of being able to implement solution $\pmb{x}^k$ might also be associated with  costs that are independent of the scenario realization. The existence of first-stage costs is a common feature in multi-stage robust and stochastic optimization. For example, in disaster management preparatory actions are necessary in order to be able to eventually deploy a specific emergency response plan. In order to reflect this option in the model we enlarge the domain to contain a variable vector $\pmb{y}\in\mathbb{R}^\pi$ that is associated to preparatory actions that have to be adequate for all potential future solutions $\pmb{x}^k$. In particular for all solutions $k$ the pair $(\pmb{y},\pmb{x}^k)$ has to part of the new domain $\bar{\mathcal{X}}\subseteq \mathbb{R}^\pi \times  \mathcal{X}$, ensuring that the performed preparations enable the realization of solution $\pmb{x}^k$. By adding the cost term $\hat{\pmb{c}}^\top\pmb{y}$ to the objective function and modeling the connection between preparation and implementation in $\bar{\mathcal{X}}$ the original model can easily be adapted to deal with such a factor.

\section{Computational Experiments}
\label{sec:exp}

We present three experiments to address the following questions:
\begin{enumerate}
\item How does the performance of the greedy heuristic compare to an optimal solution for the interpretable optimization problem?

\item How large is the loss of performance when using the \DR instead of unrestricted solution assignments in the min-sum-min setting?

\item How does our interpretable approach perform on real-data?
\end{enumerate}
All code is implemented in C++ using Cplex version 12.8 to solve mixed-integer programs. We use a compute server running Ubuntu 18.04 with ten Intel Xeon Gold 5220 CPUs running at 2.20GHz. Each process was restricted to a single thread. In all experiments, we focus on the sum of objective values as the decision criterion. Code and data are available online on GitHub\footnote{
 \url{https://github.com/goerigk/explainable-opt-code.git}}.

\subsection{Experiment 1}

To compare the performance of the optimal interpretable solution found by solving the mixed-integer programming formulation proposed in Section~\ref{sec:formulation} with the performance of the greedy heuristic from Section~\ref{sec:greedy}, we generate shortest path problems on grid graphs. We create small quadratic grids of $5\times 5$ nodes. Each node is connected to its up to four horizontal and vertical neighbors, i.e., there are 40 edges. Edges are oriented upwards and to the right. The origin is placed on the bottom left corner, while the destination is placed on the top right corner. To create scenarios, we generate fives types $t$, each given by its own midpoint $\mu^t_e$ for each edge $e$ and a deviation range $\delta^t_e$. Midpoints are sampled uniformly i.i.d. from $[10,30]$ and deviations from $[0,0.25]$. To generate a scenario, we first choose a type $t$ and then sample each edge costs uniformly i.i.d. from $[(1-\delta^t_e)\mu^t_e, (1+\delta^t_e)\mu^t_e]$. This means that each problem instance contains five types of scenarios, which enables us to test the degree to which optimization models are able to make use of this structure. We consider instances with $\{5,6,\ldots,20\}$ scenarios for training and 1000 scenarios for testing. For each number of scenarios, 100 instances are created (a total of $16\times 100 = 1600$ instances).

For each instance, we solve the optimization model~\eqref{Eq::MainModel} (denoted as IP2 and IP4 when using $K=2$ and $K=4$ solutions, respectively) and we run the greedy heuristic (denoted as G2 and G4). We compare the solution time of each approach, as well as the performance on training and test instances of the resulting \DRs. We measure performance in the following way. Let $f(\xi)$ denote the objective value of the solution found by some method in scenario $\xi$. Let $f^{nom}(\xi)$ denote the objective value of the nominal solution, i.e., the solution that minimizes the sum (or average) of objective values over all training scenarios. Finally, let $f^{opt}(\xi)$ denote the optimal objective value for scenario $\xi$. The performance of $f$ in scenario $\xi$ is then calculated as \[ 1 - \frac{ f(\xi) - f^{opt}(\xi) }{ f^{nom}(\xi) - f^{opt}(\xi) } \]
which means that a value of 100\% gives an optimal solution in scenario $\xi$, whereas a value of 0\% means that there is no improvement over the nominal solution. We consider the average of these performance values over all scenarios.

\begin{figure}[htbp]
\begin{center}
\subfigure[Training instances.\label{fig:exp1in}]{\includegraphics[width=0.48\textwidth]{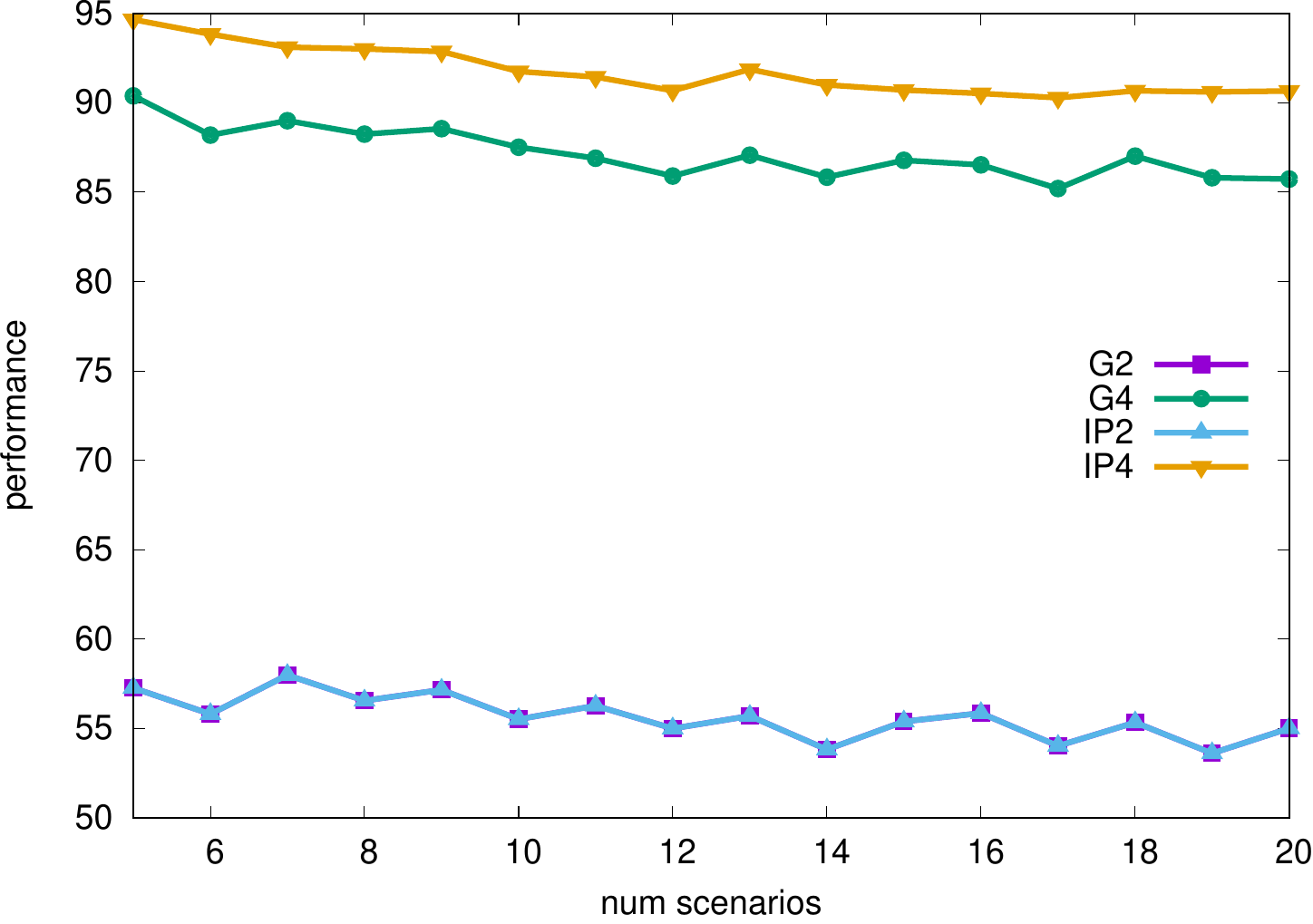}}%
\hfill
\subfigure[Test instances.\label{fig:exp1out}]{\includegraphics[width=0.48\textwidth]{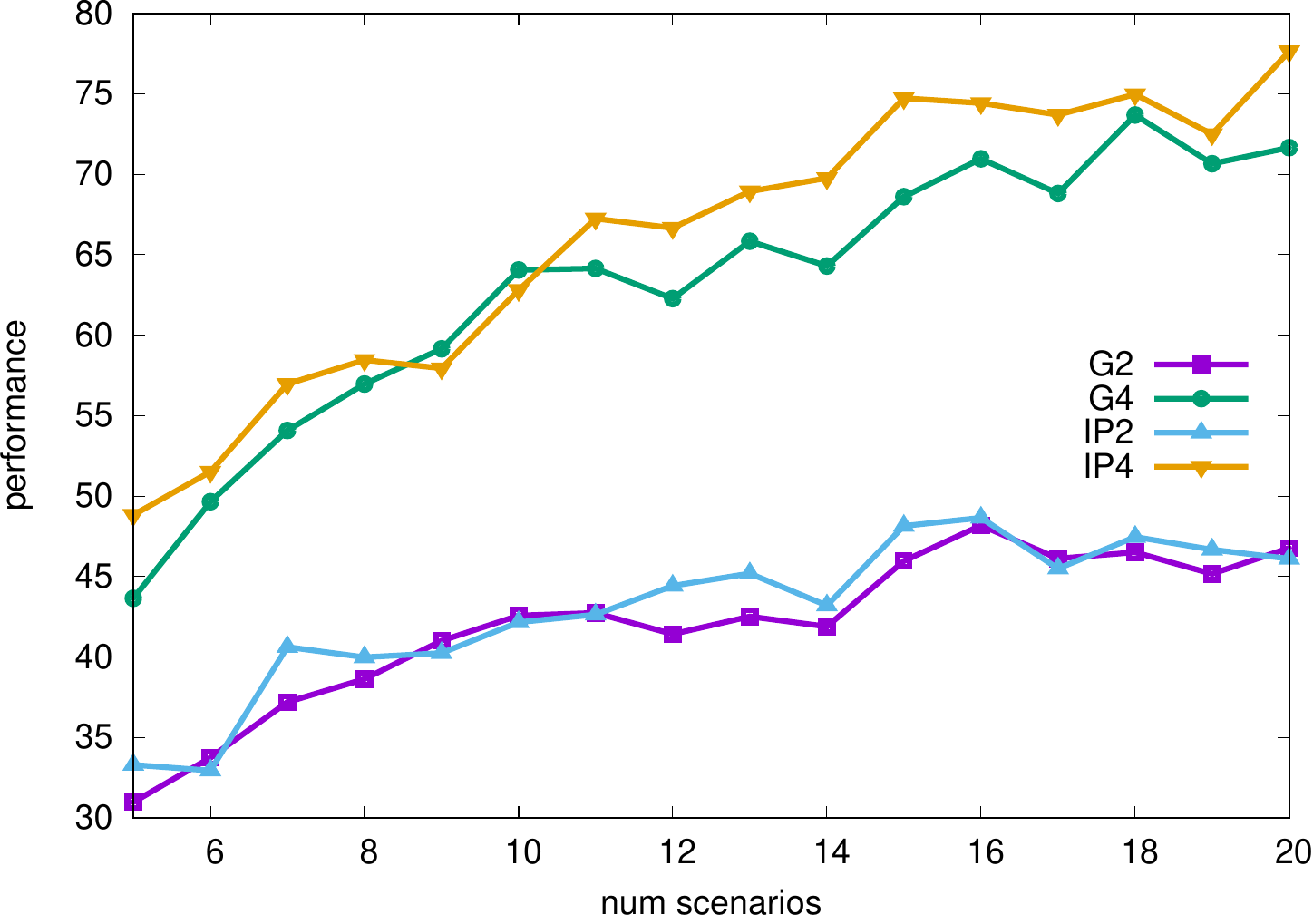}}%
\end{center}
\caption{Experiment~1, average performance values for varying numbers of training scenarios.\label{fig:exp1}}
\end{figure}

Figure~\ref{fig:exp1} shows the average performance for training scenarios and test scenarios, for different numbers of training scenarios. Note that it can be expected that the performance on training scenarios has a tendency to decrease with more training scenarios, as more scenarios need to be considered simultaneously. At the same time, more training scenarios should lead to a better performance on test scenarios. Our results confirm this behavior.

In Figure~\ref{fig:exp1in}, G2 and IP2 have identical performance, and thus both lines are on top of each other. This is to be expected, as the greedy heuristic places a single decision in an optimal way. Comparing G4 and IP4, there is a performance difference of about 5 percentage points on the training data.

On test instances (see Figure~\ref{fig:exp1out}), it is possible that G2 and IP2 have different performance values, as optimal solutions are not necessarily unique. There is no clear preference of one method over the other. We note that IP4 again outperforms G4, but the difference between both methods has become smaller. In some cases, G4 even shows a slightly better average performance than IP4. Overall, the greedy heuristic gives a reasonable approximation of the optimal solution on these instances.

\begin{figure}[htbp]
\begin{center}
\includegraphics[width=0.48\textwidth]{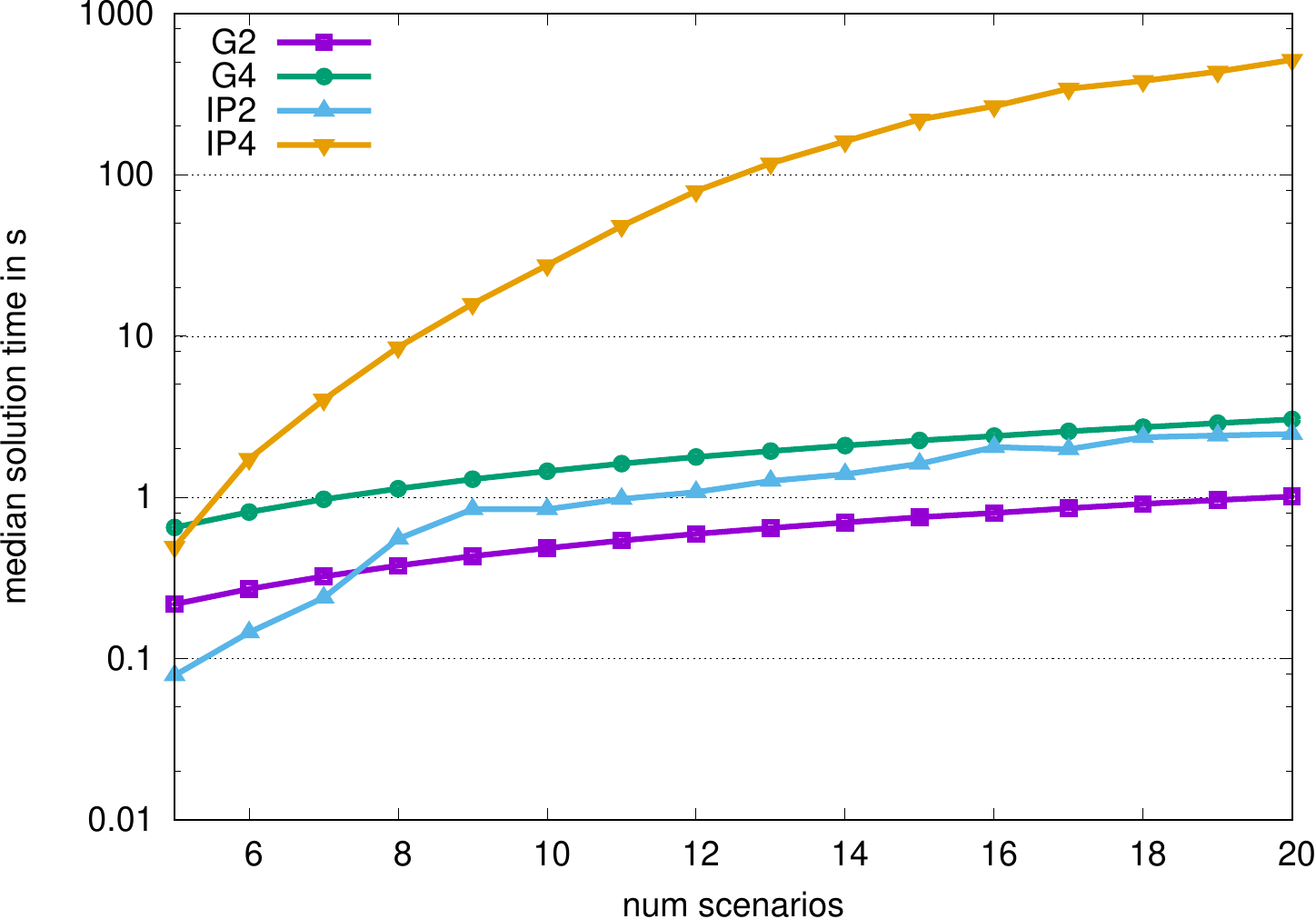}
\end{center}
\caption{Experiment~1, median computation times for varying numbers of training scenarios.\label{fig:exp1times}}
\end{figure}

Compare this relatively small loss in performance with the improvements in median computation times, as shown in Figure~\ref{fig:exp1times}. Note the logarithmic vertical axis. 

While the greedy method scales well with an increased number of scenarios, this is not the case for the integer program. For 20 scenarios (the largest problems considered in this experiment), IP4 is around two orders of magnitude slower than G4. Also note that plotting median times benefits the IP methods, as they are more prone to outliers than the greedy methods. In summary, the greedy method gives a good trade-off between small computation times and larger performance values, and will thus be the method we consider in our subsequent experiments.

\subsection{Experiment 2}

In this second experiment, we consider the possible disadvantage resulting from using \DRs over formulations that do not consider this restriction. To this end, we use the greedy heuristic to find \DRs, and compare with the performance of the min-sum-min model with two (MSM2) and with four solutions (MSM4). Performance is measured in the same way as in Experiment~1. We again use grid graphs, but consider larger instances with 10 times 10 nodes (180 edges) and a number of training scenarios in the range $\{10,15,\ldots,50\}$. We generate 100 instances for each number of training scenarios (a total of $9\times 100 = 900$ instances). For each instance, we first run the greedy heuristics G2 and G4 and then assign the computation times as time limits for MSM2 and MSM4, respectively. While this may create a bias towards the greedy methods (which have exactly as much computation time as they need), this approach avoids setting an arbitrary time limit for the MSM methods.

\begin{figure}[htbp]
\begin{center}
\subfigure[Training instances.\label{fig:exp2in}]{\includegraphics[width=0.48\textwidth]{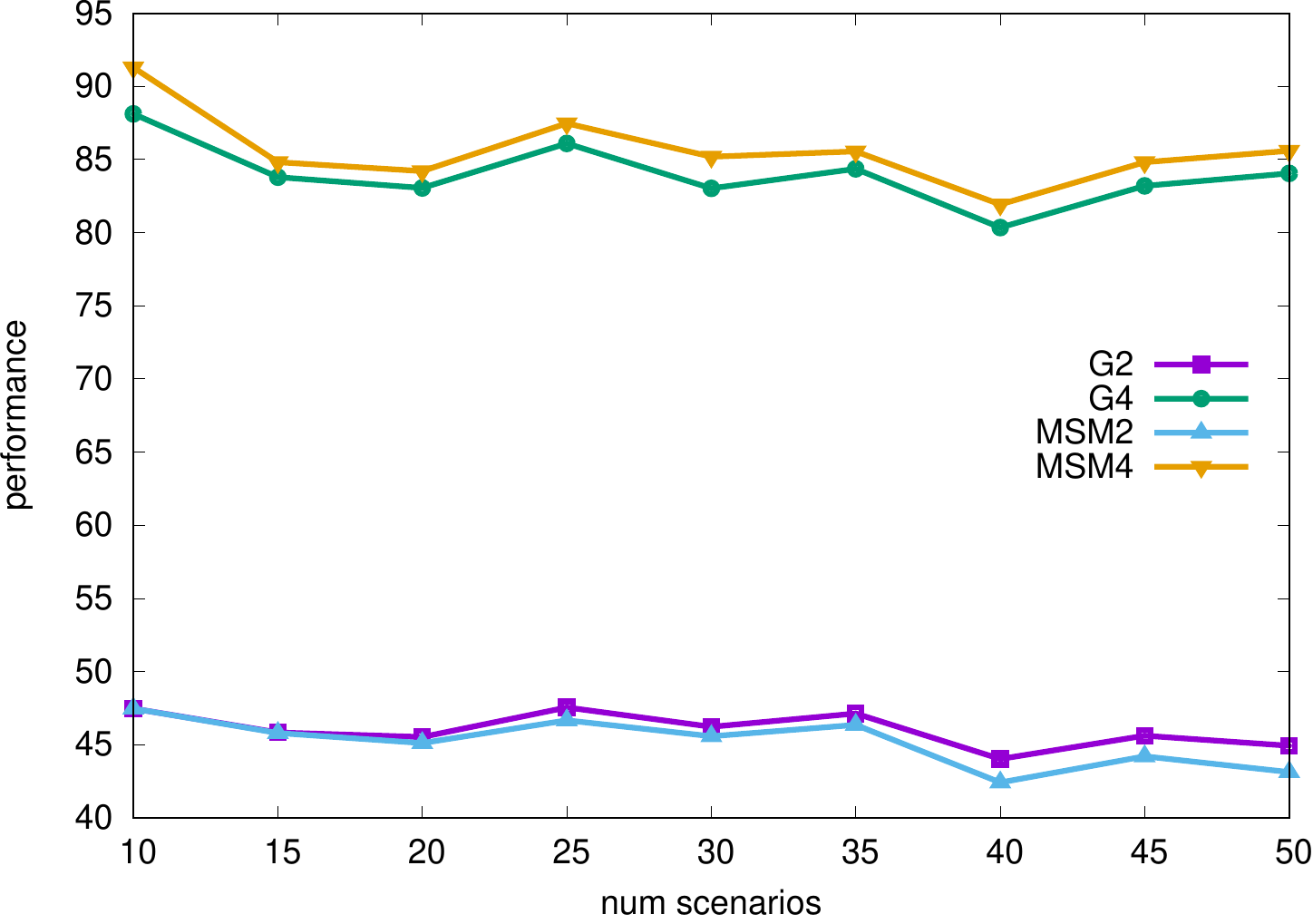}}%
\hfill
\subfigure[Test instances.\label{fig:exp2out}]{\includegraphics[width=0.48\textwidth]{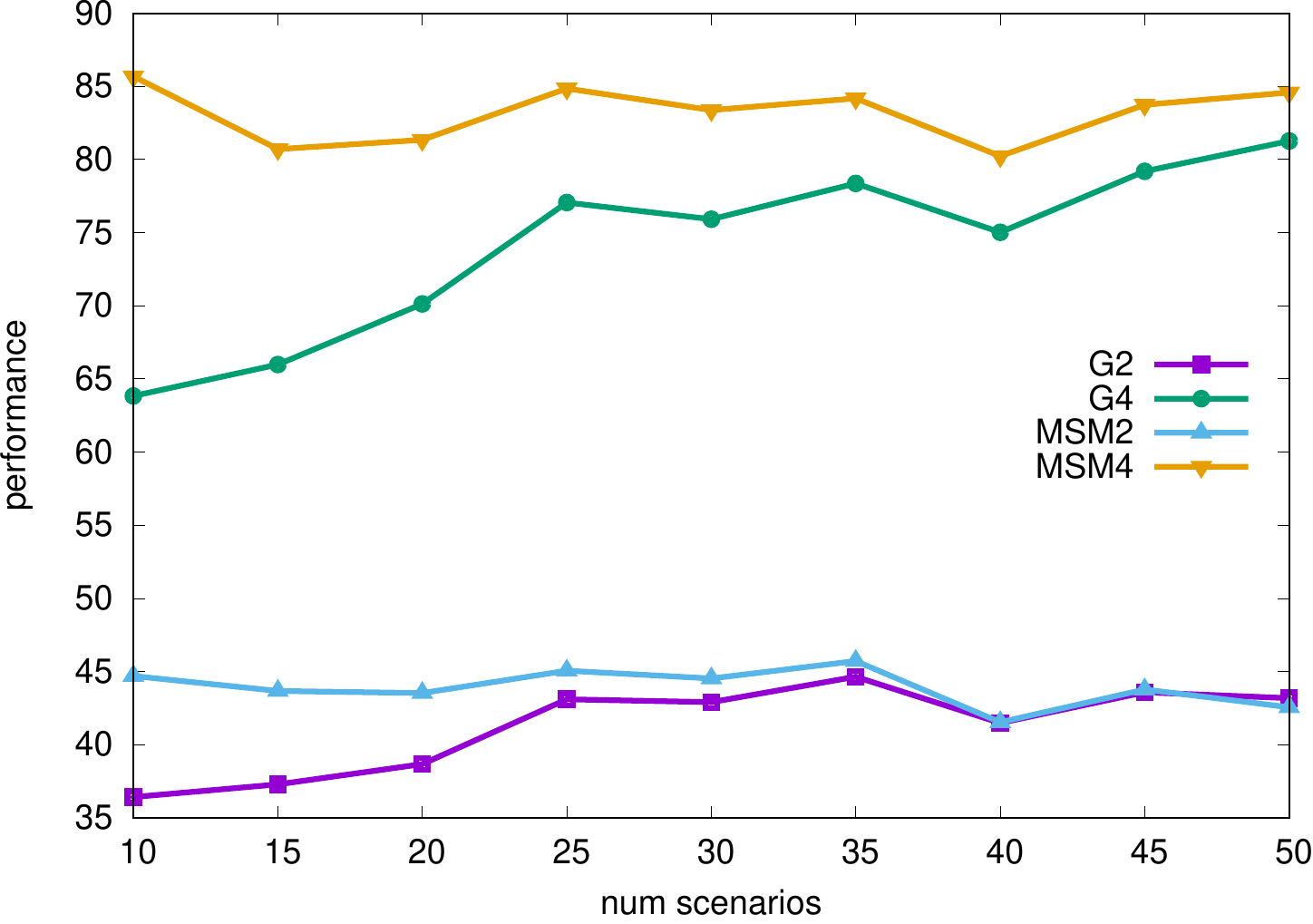}}%
\end{center}
\caption{Experiment~2, average performance values for varying numbers of training scenarios.\label{fig:exp2}}
\end{figure}

In Figure~\ref{fig:exp2}, we present the average performance on training and on test instances. G4 achieves nearly the same performance of MSM4 on training instances (see Figure~\ref{fig:exp2in}), with a performance difference of about one to two percentage points. Interestingly, G2 even outperforms MSM2, which is only possible due to the time limit applied to the latter method.

On test instances (see Figure~\ref{fig:exp2out}), we note that the usage of our \DRs show a drop in performance for small sizes of training sets. This indicates that the rules learned on these sets do not yet generalize well to other scenarios. As the number of scenarios in the training set increases, the performance of greedy again comes close to that of MSM. While G2 reaches the same performance as MSM2, the difference between G4 and MSM4 remains at around three percentage points.

\begin{figure}[htbp]
\begin{center}
\includegraphics[width=0.48\textwidth]{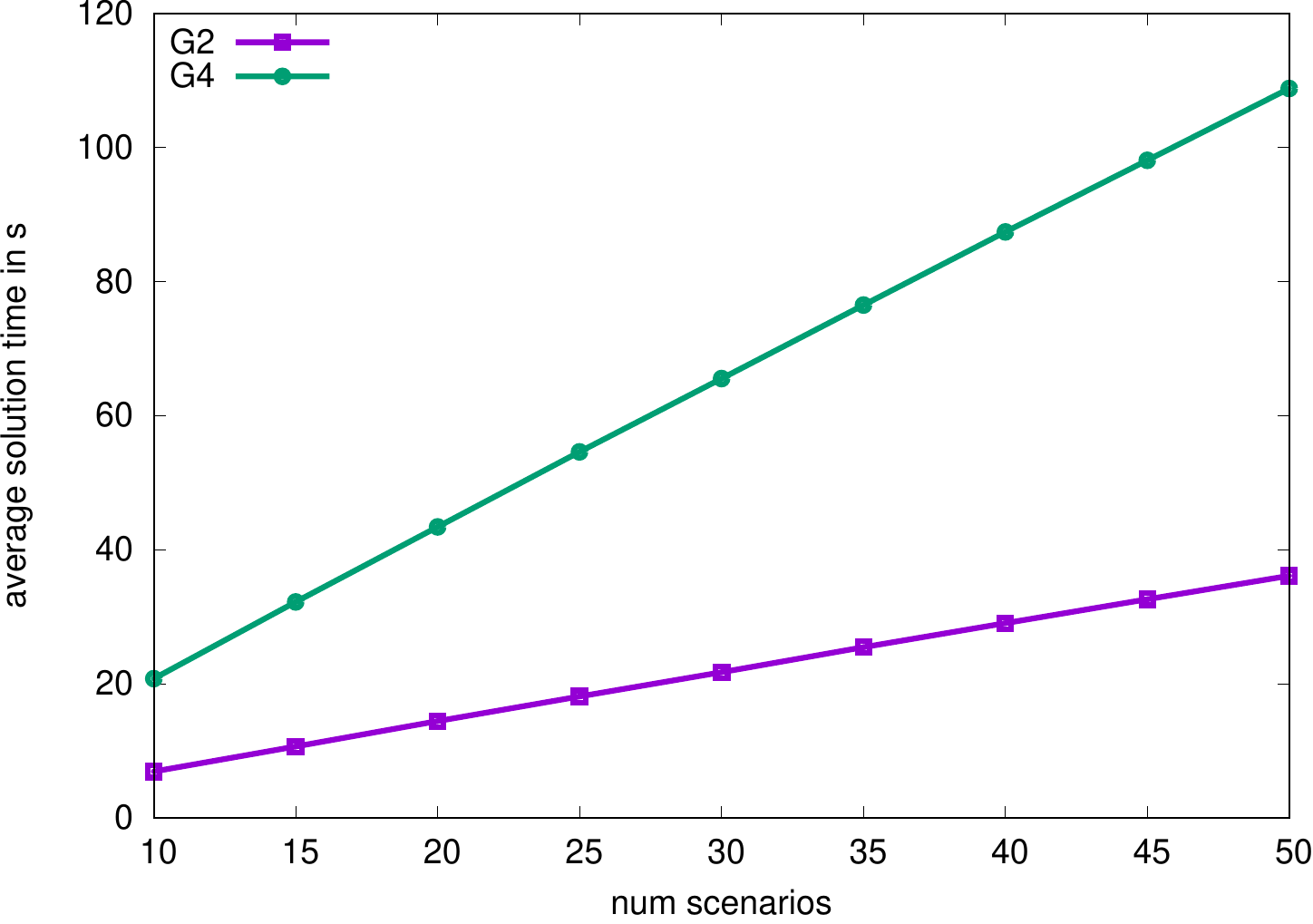}
\end{center}
\caption{Experiment~2, median computation times for varying numbers of training scenarios.\label{fig:exp2times}}
\end{figure}

We show median computation times of the greedy methods in Figure~\ref{fig:exp2times}. As predicted in our analysis, there is a clear linear relationship between the number of scenarios in the training set and computation time. Also note that running G4 takes about three times as long as running G2 (note that G4 needs to first run G2, and then run an additional iteration where four nominal problems need to be solved for every split candidate instead of two problems, which results in three times the effort of G2).

Summarizing these results, we find that \DRs only lead to a small decrease in performance compared to general decision rules on the instances we considered, in particular, if sufficient training scenarios are available. At the same time, the linear solution time of the greedy method means that our approach scales well in problem size.

\subsection{Experiment 3}

In this third experiment, we consider real-world data collected in \cite{chassein2019algorithms}. The graph models the road network of Chicago with 538 nodes and 1308 edges. Travel speeds were taken from a live traffic interface over a 46-day time horizon. The dataset we use contains all measurements that were made (a total of 4363 scenarios). These are split by using one half for training, and one half for testing. We generate 100 random source-sink pairs. To avoid path problems where both nodes are too close to each other, we generate a nominal path for each source-sink candidate, and repeat the sample while the number of edges in the nominal path is less or equal to 50.

To reduce computation times, we randomly skip checking each split candidate on the existing 1308 edges with 95\% probability, which means that the greedy heuristic is sped up by factor 20. Additionally, we extend the instance by two artificial edges to capture meta data on scenarios. These edges loop into an additional, separate node, so that they are never part of a path. On the first edge, we use the day of the week of a scenarios as its lengths, encoded as integer numbers from one (Monday) to seven (Sunday). On the second edge, we use the time within the day, encoded as seconds since midnight, i.e., an integer from 0 to 86399.

We use the greedy heuristic to find \DRs using all 1310 edges, and also using only the two edges carrying meta data. The median computation time for the first setting is 1209.5 seconds (two solutions) and 3450.5 seconds (four solutions), while the second setting results in median computation times of 126.4 seconds (two solutions) and 363.1 seconds, respectively. 

\begin{table}[htb]
\caption{Experiment 3, average performance values.\label{tab:exp3perf}}

\begin{center}
\begin{tabular}{r|rr|rr}
& \multicolumn{2}{|c|}{all edges} & \multicolumn{2}{|c}{meta only} \\
Dataset  & G2 & G4 & G2 & G4 \\
\hline
Training & 14.0\% & 20.8\%  & 2.8\% & 3.8\% \\
Test & 12.1\% & 17.6\% & 1.8\% & 2.7\%
\end{tabular}
\end{center}
\end{table}

\begin{figure}[htbp]
\begin{center}
\includegraphics[width=0.6\textwidth]{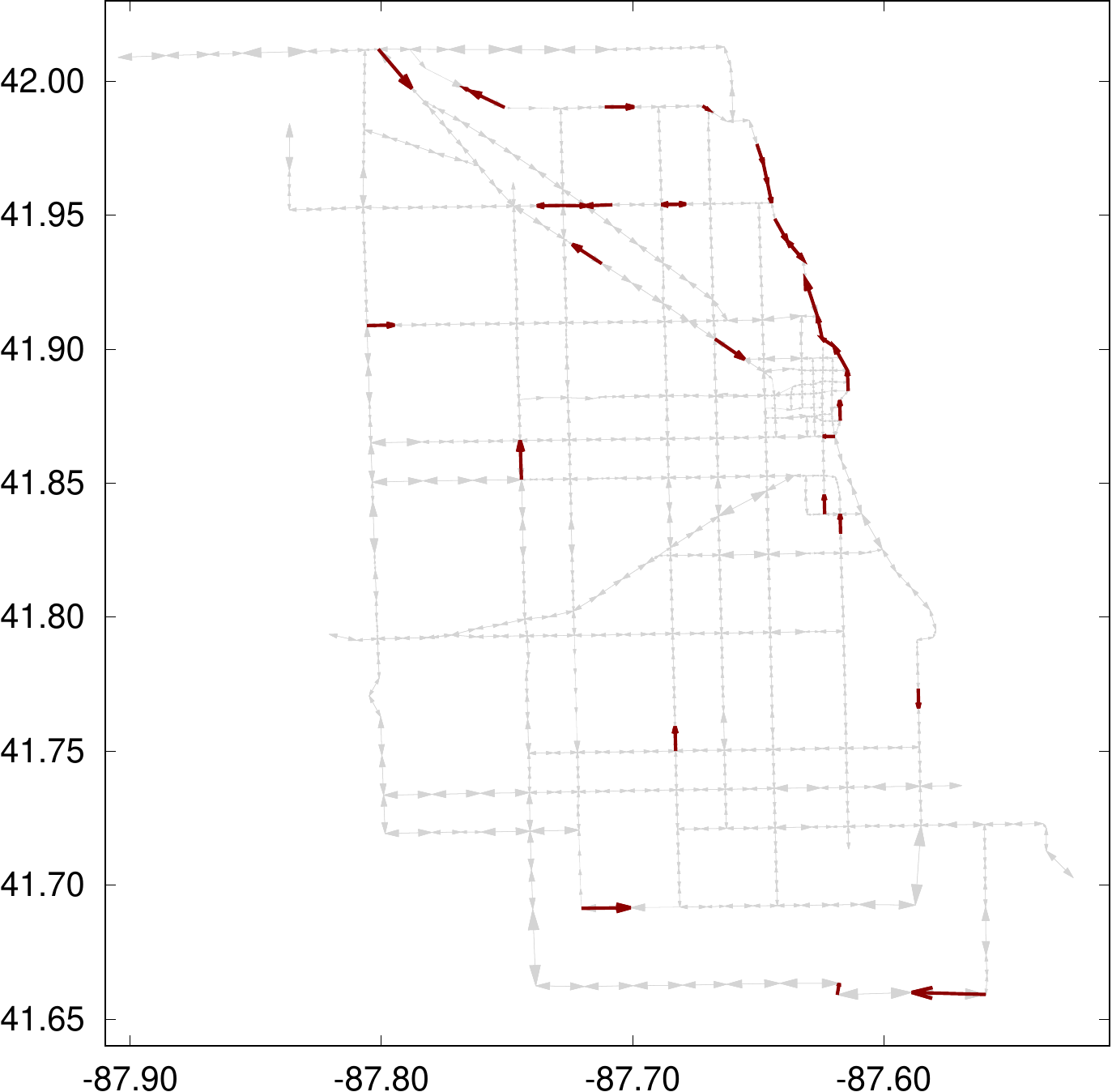}
\caption{Experiment 3, edges that were used in \DRs.}\label{fig:exp3}
\end{center}
\end{figure}

In Table~\ref{tab:exp3perf}, we show the average performance for both settings. While it is possible to slightly improve upon the nominal solution using only meta data, this improvement is small. Significantly larger improvements are possible when using the full data of the graph. Indeed, if we consider the first splitting decisions of this approach (i.e., the rule that G2 finds), only 1 in 100 runs use the day of the week (split between Monday to Saturday and Sunday), and 8 in 100 runs use the time of the day. In Figure~\ref{fig:exp3}, we highlight which edges were split on in the remaining 91 cases. We note that in particular the edges along the shore of Lake Michigan in the north-east are used for this purpose. A possible explanation is that these edges see continuous high-density traffic, which makes them reliable predictors to identify the traffic scenario.

\begin{figure}[htbp]
\begin{center}
\includegraphics[width=0.6\textwidth]{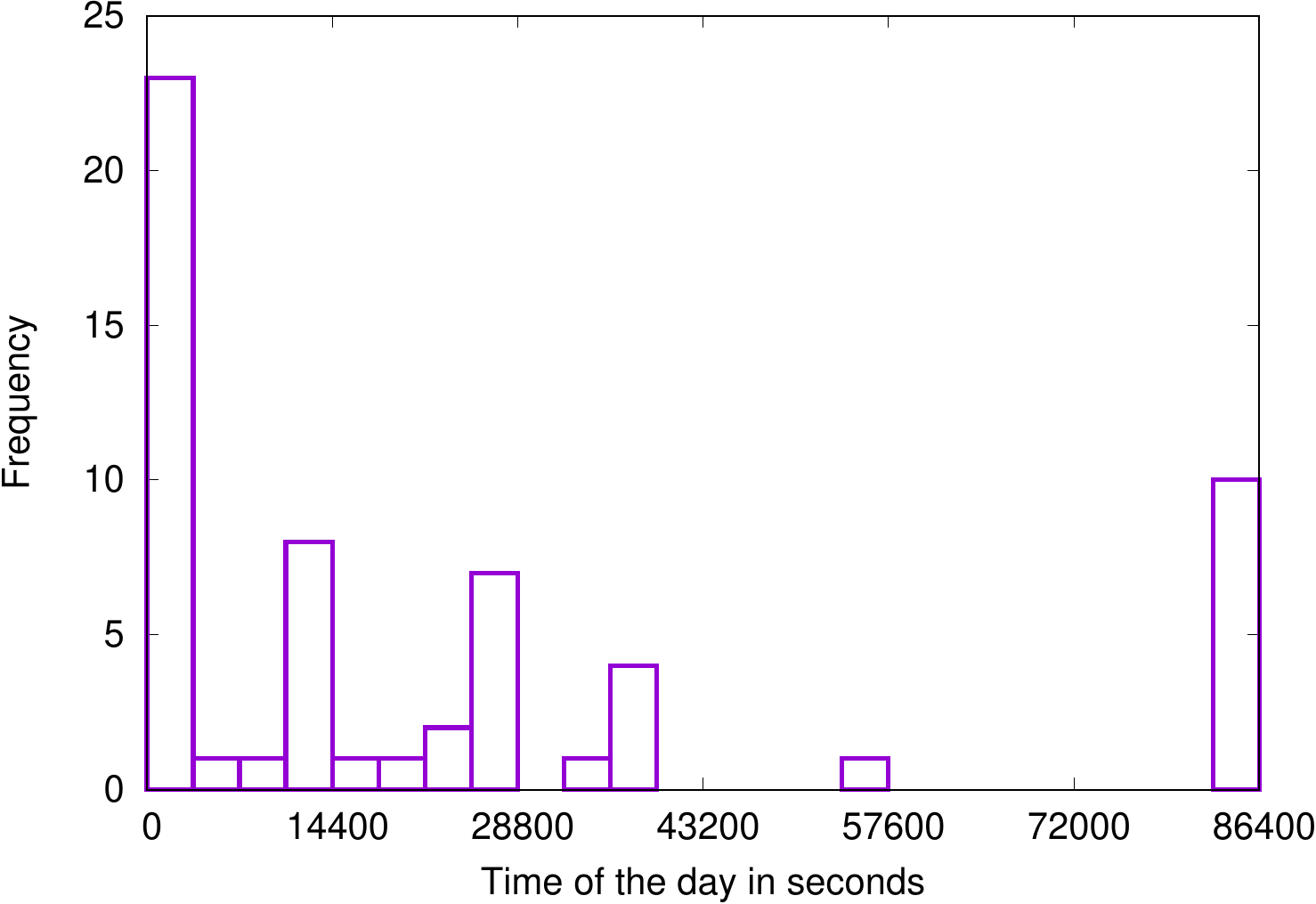}
\caption{Experiment 3, histogram of threshold values for splits on the time dimension.}\label{fig:exp3time}
\end{center}
\end{figure}

If we only consider the meta data, there are five cases where G2 splits on the day (once between Friday and Saturday, four times between Saturday and Sunday). The remaining splits are made on the time of the day. In Figure~\ref{fig:exp3time}, we present a histogram of these decisions, where each bin represents an hour of the day. The majority of splits cut off few scenarios on the boundary of the day, in particular in the first and in the last hour. This underlines the small improvements that these rules give, see Table~\ref{tab:exp3perf}. Traffic speed is highly volatile, and the data represents snapshots at particular points in time, rather than averages over larger time horizons. Hence, \IDRs such as ''before 7am'' and ''after 7am'' do not perform as well as one might intuitively expect.

\subsection{Discussion}

Our computational experiments have been designed to answer three questions. The first is concerned with a comparison of the greedy heuristic and the mixed-integer programming formulation to find \DRs in form of a decision tree. Recent research on decision trees for classification models (see \cite{bertsimas2017optimal}) reports promising results when constructing optimal decision trees. In our setting, computation times quickly increase to become prohibitive for larger problems. The greedy heuristic scales well (linearly in the number of scenarios) and is a viable choice also for larger optimization problems. At the same time, the quality of solutions found by the heuristic is close to that of exact solutions, which makes the heuristic the preferred solution method for the remaining experiments.

The second question considers the loss that occurs when we use an \DR instead of an arbitrary rule to make decisions. Our experiments show that on training data, there is little difference between both approaches. This is different on test data, where an insufficient training size leads to interpretable solutions that do not perform well. With more training data available, we can again observe that \DRs show a performance that is very close to more general rules. Hence, these results seem to encourage the use of \IDRs.

Finally, we consider the performance of our approach on real-world data. While the previous experiments constructed scenarios is a way that there exists inherent structure that can be exploited, this is not necessarily the case for the considered data set. Indeed, improvements over nominal solutions towards the lower bound are significantly smaller, but still useful with up to around 20\%. It seems that particularly high-density edges give reliable predictions which of the candidate solutions to use.

The results presented here are, naturally, only valid within the scope of the setup we used. At the same time, none of our methods were specifically designed for the shortest path problem, which underlies all experiments. Hence, we can expect to find a similar performance of our approach for different optimization problems.

\section{Conclusions}

While immense progress has been made towards finding optimal solutions to decision making problems efficiently, even the best solution is only useful if it can be applied in practice and is accepted by the stakeholders. A major barrier in this process is the black-box nature of optimization solvers. Interpretability and explainability, concepts that have become central in artificial intelligence, have seen significantly less attention in operations research and management science.

In this paper, we propose an inherently interpretable optimization model to derive solutions that come with an \DR, under which circumstances which solution should be implemented. We focused on decision tree based rules, which are easily comprehensible and applicable.

To find such solutions and \IDRs, a mixed-integer programming formulation can be used. For decision trees with small depth, it is also viable to enumerate candidate solutions. Finally, it is possible to construct heuristic rules greedily. This basic model can be easily extended to include a range of further features.

Our experiments indicate that the greedy heuristic can be used to find high-quality solution in reasonable time. Moreover, it is possible to include any available nominal solution algorithm in the process, which improves the applicability of our method. We found that the loss of using \IDRs instead of complex, but incomprehensible rules is in the range of only a few percentage points.

As few methods for interpretable optimization exist, there are many avenues for further research that should be explored. In this paper, we focussed on decision trees as \IDRs. Other models can be used instead, such as point-based systems, where a decision maker only needs to calculate a score by summing up points for each feature that is checked (i.e., using a preferably sparse perceptron). Furthermore, these methods should be tested in practice with decision makers to see which rules find acceptance. Finally, we considered a setting where a set of candidate solutions is found to cover any potential optimization instance that we may encounter. Fundamentally different approaches are conceivable, for example, we might be interested to produce just one solution for one particular problem, but this solution is required to have a simple and communicable structure. Imagine a car driver who asks for directions. He or she would be interested not necessarily in the fastest path to reach the destination, but rather in a path that is easy to find, i.e., a path with few turns that mostly uses major roads.

\end{document}